\documentclass[a4paper,12pt]{article}

\usepackage{amsthm}
\usepackage{amsmath,amssymb,latexsym,amsfonts,mathrsfs}

\renewcommand{\d}{{\rm d}} %differential
\newcommand{\e}{{\rm e}} 
\newcommand{\dist}{\stackrel{{\rm d}}{=}}
\newcommand{\cdist}{\stackrel{{\rm d}}{\longrightarrow}}
\newcommand{\tend}[2]{\mathrel{\mathop{\longrightarrow}\limits^{#1}_{#2}}}
\renewcommand{\hat}{\widehat}
\renewcommand{\tilde}{\widetilde}

\newcommand{\up}{\ensuremath{\uparrow}}
\newcommand{\down}{\ensuremath{\downarrow}}

\newcommand{\rbra}[1]{\!\left( #1 \right)} %round brackets or parentheses
\newcommand{\cbra}[1]{\!\left\{ #1 \right\}} %curly brackets or braces
\newcommand{\bE}{\ensuremath{\mathbb{E}}}
\newcommand{\bP}{\ensuremath{\mathbb{P}}}
\newcommand{\bR}{\ensuremath{\mathbb{R}}}
\newcommand{\bZ}{\ensuremath{\mathbb{Z}}}

\newcommand{\cA}{\ensuremath{\mathcal{A}}}
\newcommand{\cB}{\ensuremath{\mathcal{B}}}

\newcommand{\cG}{\ensuremath{\mathcal{G}}}

\newcommand{\cN}{\ensuremath{\mathcal{N}}}

\newcommand{\vs}{\ensuremath{\boldsymbol{s}}}

\theoremstyle{plain}
\newtheorem{Thm}{Theorem}[section]
\newtheorem{Lem}[Thm]{Lemma}
\newtheorem{Prop}[Thm]{Proposition}

\theoremstyle{definition}

\newcommand{\Proof}[2][Proof]{\begin{proof}[{#1}] #2 \end{proof}}

\numberwithin{equation}{section}

\makeatletter
\renewcommand\section{\@startsection {section}{1}{\z@}%
                                   {-3.5ex \@plus -1ex \@minus -.2ex}%
                                   {2.3ex \@plus.2ex}%
                                   {\normalfont\large\bf}}
\makeatother

\makeatletter
\renewcommand\subsection{\@startsection {subsection}{1}{\z@}%
                                   {-3.5ex \@plus -1ex \@minus -.2ex}%
                                   {2.3ex \@plus.2ex}%
                                   {\normalfont\normalsize\bf}}
\makeatother

\usepackage[dvipdfm]{graphicx}
\usepackage{float}

\usepackage[usenames]{color}

\begin{document}

\begin{center}
{\Large \bf 
Arcsine and Darling--Kac laws for piecewise linear random interval maps 
}
\end{center}
\begin{center}
Genji \textsc{Hata}\footnote{
SAPIX 
} \quad and \quad 
Kouji \textsc{Yano}\footnote{
Graduate School of Science, Kyoto University.}\footnote{
The research of Kouji Yano was supported by 
JSPS KAKENHI grant no.'s JP19H01791 and JP19K21834 
and by JSPS Open Partnership Joint Research Projects grant no. JPJSBP120209921. 
This research was supported by RIMS and by ISM.}
\end{center}

\begin{abstract}
We give examples of piecewise linear random interval maps 
satisfying arcsine and Darling--Kac laws, 
which are analogous to Thaler's arcsine and Aaronson's Darling--Kac laws 
for the Boole transformation. 
They are constructed by random switch of two piecewise linear maps 
with attracting or repelling fixed points, 
which behave as if they were indifferent fixed points of a deterministic map. 
\end{abstract}


\section{Introduction} \label{sec: intro}

Let $ \cA $ and $ \cN $ denote arcsine and standard normal random variables, i.e., 
\begin{align}
\bP(\cA \in \d u) = \frac{\d u}{\pi \sqrt{u(1-u)}} \ (0<u<1) 
, \quad 
\bP(\cN \in \d u) = \e^{-u^2/2} \frac{\d u}{\sqrt{2 \pi}} \ (u \in \bR) . 
\label{}
\end{align}
It is well-known that a simple symmetric random walk 
$ \cbra{ W_n }_{n=0}^{\infty } $ on $ \bZ $ 
satisfies \emph{L\'evy's arcsine law} and \emph{Darling--Kac law} 
(\cite{MR0000919}; see also \cite{MR84222}), respectively: 
\begin{align}
\frac{1}{N} \sum_{n=0}^{N-1} 1_{\{ W_n>0 \}} 
\tend{\rm d}{N \to \infty} \cA 
\quad \text{and} \quad 
\frac{1}{\sqrt{N}} \sum_{n=0}^{N-1} 1_{\{ W_n \in E \}} 
\tend{\rm d}{N \to \infty} \frac{\# E}{\sqrt{\pi}} \, |\cN| 
\label{}
\end{align}
for all bounded set $ E \subset \bZ $, where $ \# E $ denotes the number of elements of $ E $ 
and $ \cdist $ means the convergence in distribution.

Thaler and Aaronson 
obtained analogous results for the \emph{Boole transformation} 
$ T $ of $ [0,1] $ defined as 
\begin{align}
T(x) = \frac{x(1-x)}{1-x-x^2} \ (0<x<1/2) 
, \quad 
T(x) = 1-T(1-x) \ (1/2<x<1) . 
\label{}
\end{align}
Thaler's arcsine law (\cite{Th}) and Aaronson's Darling--Kac law (\cite[Theorem 1]{MR632462}) 
can be stated as follows: 
For any random initial point $ \Theta $ in $ [0,1] $ with a.c. density, 
it holds that 
\begin{align}
\frac{1}{N} \sum_{n=0}^{N-1} 1_{ \{ T^n(\Theta) > 1/2 \} } 
\tend{\rm d}{N \to \infty } \cA 
\quad \text{and} \quad 
\frac{1}{\sqrt{N}} \sum_{n=0}^{N-1} 1_{\{ T^n(\Theta) \in E \}} 
\tend{\rm d}{N \to \infty} \frac{\mu(E)}{\sqrt{\pi}} \, |\cN| 
\label{eq: Thaler Aaronson}
\end{align}
for all Borel set $ E $ with $ \mu(E) < \infty $, 
where $ \mu $ is the unique (up to a constant multiple) 
a.c. $ \sigma $-finite $ T $-invariant measure given as 
\begin{align}
\mu(\d x) = \Phi'(x) \d x = \rbra{ \frac{1}{x^2} + \frac{1}{(1-x)^2} } \d x 
\quad \text{on $ [0,1] $}. 
\label{}
\end{align}
We note that $\mu$ has infinite mass near 0 and 1 
and has finite mass away from 0 and 1. 
The points $ 0 $ and $ 1 $ are \emph{indifferent fixed points} for $ T $ in the sense that 
\begin{align}
[T(0+) = 0 , \ T'(0+)=1 ] \quad \text{and} \quad [ T(1-) = 1 , \ T'(1-) = 1 ]. 
\label{}
\end{align}
These facts show that the orbit of $ T $ spends most of time 
near the indifferent fixed points 0 and 1, 
and the limit laws \eqref{eq: Thaler Aaronson} 
characterize the ratio of time spent near 1 
and the decay of the ratio of time spend away from 0 and 1. 
(Note that our transformation $ T $ of $ [0,1] $ 
can be obtained, via the change of variables $ y = \frac{2x-1}{x(1-x)} $, 
from the original Boole transformation of $ \bR $ defined as $ S(y) = y - 1/y $ 
(see \cite{MR335751}), 
which preserves the Lebesgue measure on $ \bR $.)

Our aim is to obtain arcsine and Darling--Kac laws for random maps 
analogous to \eqref{eq: Thaler Aaronson}. 
For two deterministic interval maps $ \tau_1, \tau_2 : [0,1] \to [0,1] $ 
and a constant $ 0 < p < 1 $, 
we consider the random map 
\begin{align}
T = 
\begin{cases}
\tau_1 & (\text{with probability $ p $}) , 
\\
\tau_2 & (\text{with probability $ 1-p $}) . 
\end{cases}
\label{eq: randommap}
\end{align}
A measure $ \mu $ on $ [0,1] $ is called \emph{$ T $-invariant} 
if $ \mu $ is not a zero measure and
\begin{align}
(\ \bE \mu \circ T^{-1} = \ ) \ p \mu \circ \tau_1^{-1} + (1-p) \mu \circ \tau_2^{-1} = \mu . 
\label{eq: T-inv}
\end{align}
Let $ \cbra{ T_n }_{n=1}^{\infty } $ be an i.i.d. sequence of random maps 
with $ T_1 \dist T $ 
and we define
\begin{align}
T^{(n)} = T_n \circ T_{n-1} \circ \cdots \circ T_1. 
\label{eq: random map iteration}
\end{align}
The resulting random map $T^{(n)}$ can be regarded as 
the $ n $-fold composition of $ T $.

Based on Hata \cite{H}, 
we adopt the two deterministic maps 
\begin{align}
\tau_1(x)= \begin{cases}
x/2 & (0 < x < 1/2) \\
2x-1 & (1/2 < x < 1)
\end{cases}
, \quad 
\tau_2(x)= \begin{cases}
2x & (0 < x < 1/2) \\
(x+1)/2 & (1/2 < x < 1) 
\end{cases}
. 
\label{eq: hatamap}
\end{align}
Note that 0 for $ \tau_1 $ and 1 for $ \tau_2 $ are attracting fixed points: 
\begin{align}
\tau_1^n(x) \tend{}{n \to \infty } 0 
, \quad 
\tau_2^n(x) \tend{}{n \to \infty } 1 
\quad \text{for $ 0<x<1 $ except for $ 1/2 $}, 
\label{}
\end{align}
while 1 for $ \tau_1 $ and 0 for $ \tau_2 $ are repelling fixed points. 
We call the corresponding random map $ T $ the \emph{Hata map}. 
\begin{figure}[h]
\begin{minipage}{0.5\textwidth}
\centering
 \includegraphics[width=30ex]{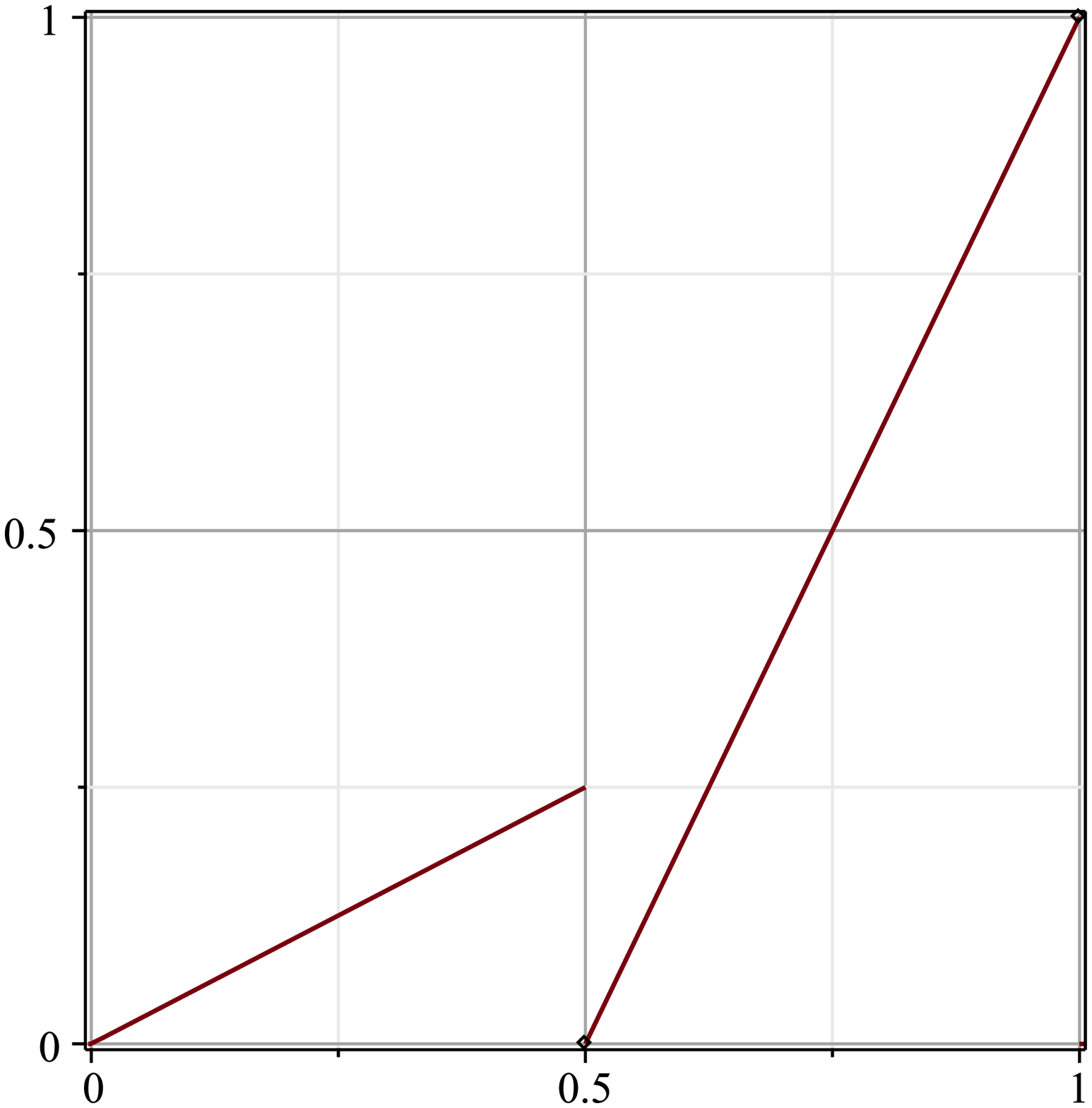}
 \caption{$\tau_1(x)$}
 \label{fig:tau1}
\end{minipage}
\begin{minipage}{0.5\textwidth}
\centering
 \includegraphics[width=30ex]{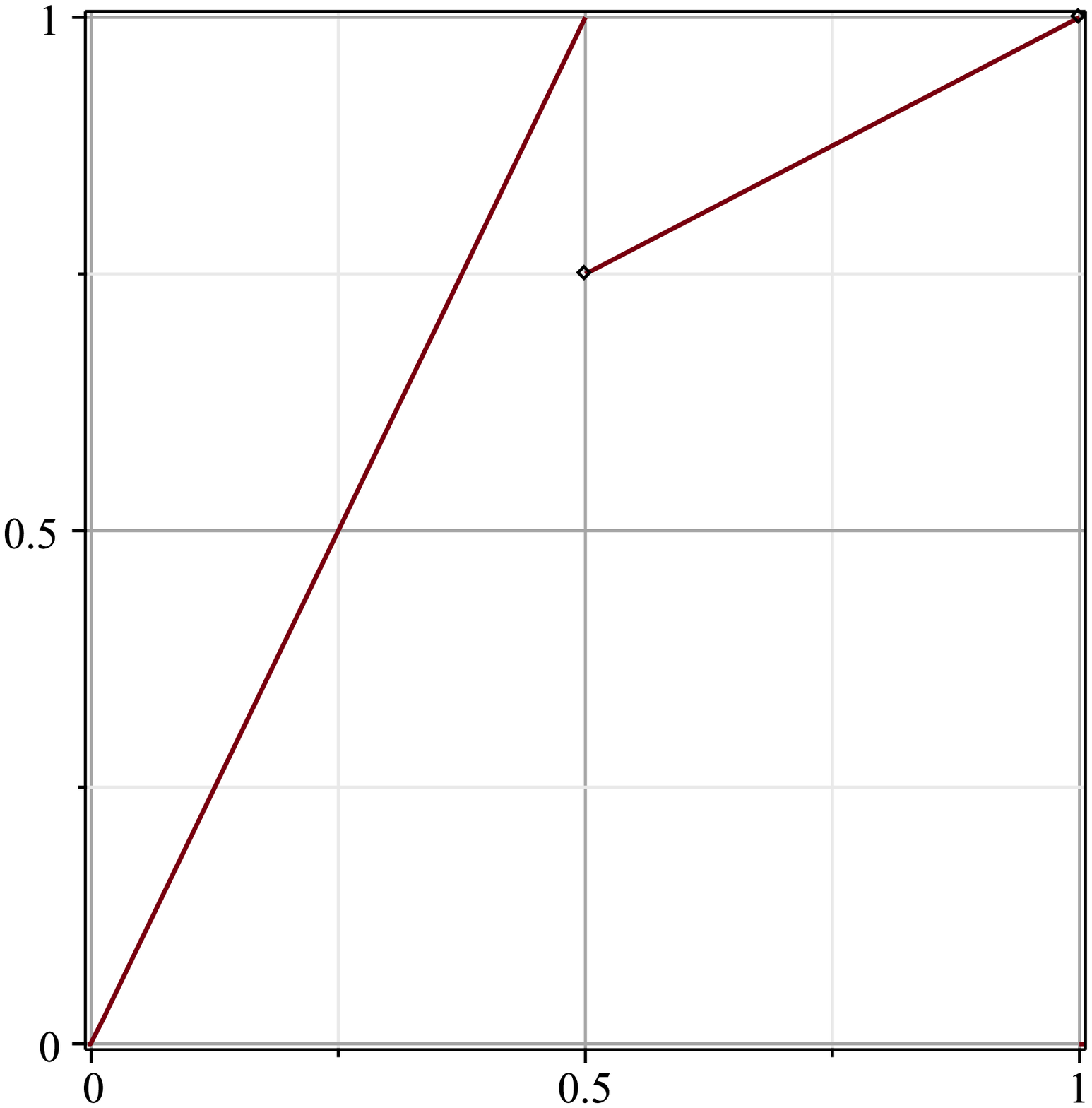}
 \caption{$\tau_2(x)$}
 \label{fig:tau2}
\end{minipage}
\end{figure}
We choose $ p = 1/2 $ so that 0 and 1 are \emph{indifferent-in-average fixed points} 
in the sense that 
\begin{align}
[T(0+)=0 , \ \bE \log |T'(0+)| = 0 ] 
\quad \text{and} \quad 
[T(1-)=1 , \ \bE \log |T'(1-)| = 0 ] . 
\label{}
\end{align}
We expect that the attracting and repelling effects are balanced in this case. 

Let $ \lambda $ denote the Lebesgue measure. 
We define $ \alpha = \cbra{ I^-_k, I^+_k }_{k=0}^{\infty } $ by 
\begin{align}
I^-_k = \rbra{ \frac{1}{2^{k+2}}, \frac{1}{2^{k+1}} } 
, \quad 
I^+_k = \rbra{ 1-\frac{1}{2^{k+1}}, 1-\frac{1}{2^{k+2}} }. 
\label{eq: alpha}
\end{align}
For the Hata map the family $ \alpha $ is a \emph{$ \lambda $-partition} of $ [0,1] $ 
in the sense that 
$ \alpha $ is a disjoint family of Borel subsets of $ [0,1] $ such that 
$ \sum_{s \in \alpha } s = [0,1] $ mod $ \lambda $ 
and $ 0 < \lambda(s) < \infty $ for all $ s \in \alpha $. 
Note that we sometimes write $ \sum $ instead of $ \bigcup $ for disjoint union. 
We say that a measure $ \mu $ is \emph{locally constant on $ \alpha $} 
if it satisfies 
\begin{align}
\mu(\d x) = \sum_{s \in \alpha } \frac{\mu(s)}{\lambda(s)} 1_s(x) \d x . 
\label{eq: im}
\end{align}

\begin{Thm} \label{arclaw}
Regarding the Hata map with $ p=1/2 $, there exists 
a unique (up to a constant multiple) 
$ \lambda $-a.c. $ \sigma $-finite $ T $-invariant measure, 
which is locally constant on $ \alpha $ with 
\begin{align}
\mu(I^-_k) = \mu(I^+_k) 
= \rbra{ 2^{k+1} - 1 } \cdot 2^{-k-2} 
, \quad k=0,1,2,\ldots 
\label{eq: im2}
\end{align}
Consequently, $ \mu $ has infinite mass near 0 and 1 
and has finite mass away from 0 and 1. 
Moreover, for any random initial point $\Theta$ in $ [0,1] $ with $ \lambda $-a.c. density 
which is independent of the random maps $ \{ T_n \} $, 
it holds that 
\begin{align}
\frac{1}{N} \sum_{n=0}^{N-1} 1_{ \{ T^{(n)}(\Theta) > 1/2 \} } 
\tend{\rm d}{N \to \infty } \cA , 
\quad \text{and} \quad 
\frac{1}{\sqrt{N}} \sum_{n=0}^{N-1} 1_{\{ T^{(n)}(\Theta) \in E \}} 
\tend{\rm d}{N \to \infty} \frac{\mu(E)}{\sqrt{\pi}} \, |\cN| 
\label{eq: arclaw}
\end{align}
for all Borel set $ E $ with $ \mu(E) < \infty $, 
where by $ \cdist $ 
we mean the convergence in distribution on the extended probability space. 
\end{Thm}

This result shows that the indifferent-in-average fixed points 0 and 1 for the random map $ T $ 
behave as if they were indifferent fixed points of a deterministic map. 
(Note that both results of \eqref{eq: arclaw} are annealed ones.) 
The proof of Theorem \ref{arclaw} will be given in Section \ref{hata}, 
which will be divided into the following steps: 
\begin{enumerate}

\item 
We show irreducible recurrence of the Markov chain on the $ \lambda $-partition. 

\item 
We modify the $ \lambda $-partition to be a Markov partition for the skew-product. 

\item 
We verify conjugacy between the random dynamical system and the Markov chain. 

\item 
We study the wandering rate asymptotics 
and check that 
the assumptions of Thaler--Zweim\"uller's theorems \cite{TZ} are satisfied. 

\end{enumerate}

While in many cases it is easy to obtain an invariant measure by reducing the problem 
to a Markov chain, there may be a technical difficulty in obtaining a Markov partition. 
Unfortunately, we cannot say anything about what happens when $ p \neq 1/2 $, 
as we cannot find a Markov partition for the skew-product. 

L\'evy's arcsine and the Darling--Kac laws 
have been extended for dynamical systems 
by Thaler \cite{Th} and by Aaronson \cite{MR632462}, respectively. 
Their results were generalized in a common framework 
by Thaler--Zweim\"uller \cite{TZ} and Zweim\"uller \cite{MR2373265}. 
Aaronson's Darling--Kac law was generalized to convergence on a functional space 
by Aaronson \cite{MR841603} and Owada--Samorodnitsky \cite{MR3298473}. 
In the case of several indifferent fixed points, 
a joint-distributional generalization of Thaler's arcsine law 
was obtained by Sera--Yano \cite{MR3988607}. 
Recently Sera \cite{MR4063962} obtained 
a functional and joint-distributional generalization 
of Thaler's arcsine and Aaronson's Darling--Kac laws.

Among others, we would like to focus on the two random maps 
with infinite invariant measures obtained 
by Pelikan \cite{MR722776} and by Boyarsky--G\'ora--Islam \cite{BGI}. 
For other results on random maps with infinite $ T $-invariant measures, 
see Bahsoun--Bose--Duan \cite{MR3225871,MR3345162}, 
Bahsoun--Bose \cite{MR3476513}, 
Gharaei--Homburg \cite{MR3600645}, 
Abbasi--Gharaei--Homburg \cite{MR3826118}, 
Inoue \cite{MR4039772}, 
Toyokawa \cite{MR4280946} 
and Homburg--Kalle--Ruziboev--Verbitskiy--Zeegers \cite{HKRVZ}. 

First, we take up 
\begin{align}
\tau_1(x) = 2x \ \text{mod} \ 1 
, \quad 
\tau_2(x) = x/2 . 
\label{}
\end{align}
and we call the corresponding random map $ T $ the \emph{Pelikan map}. 
Note that both 0 and 1 for $ \tau_1 $ are repelling fixed points 
and 0 for $ \tau_2 $ is an attracting fixed point. 
Although Pelikan \cite{MR722776} excluded the case $ p=1/2 $, 
we choose $ p = 1/2 $ so that 0 is an indifferent-in-average fixed point 
in the sense that 
\begin{align}
[T(0+)=0 , \ \bE \log |T'(0+)| = 0 ] . 
\label{}
\end{align}
We introduce the $ \lambda $-partition for the Pelikan map as 
$ \alpha = \{ I_k \}_{k=0}^{\infty } $ with 
\begin{align}
I_k = \rbra{ \frac{1}{2^{k+1}},\frac{1}{2^k} } . 
\label{eq: P alpha}
\end{align}
Although the arcsine law does not make sense, 
we obtain the Darling--Kac law as follows. 

\begin{Thm} \label{P DKlaw}
Regarding the Pelikan map with $ p=1/2 $, there exists 
a unique (up to a constant multiple) $ \lambda $-a.c. $ \sigma $-finite $ T $-invariant measure, 
which is locally constant with 
\begin{align}
\mu(I_k) = 2^{-2} \cdot \rbra{ 2^{k+1} - 1 } \cdot 2^{-k-1} 
, \quad k=0,1,2,\ldots 
\label{eq: mm}
\end{align}
Consequently, $ \mu $ has infinite mass near 0 
and has finite mass away from 0. 
Moreover, for any random initial point $\Theta$ in $ [0,1] $ with $ \lambda $-a.c. density 
which is independent of the random maps $ \{ T_n \} $, 
it holds that 
\begin{align}
\frac{1}{\sqrt{N}} \sum_{n=0}^{N-1} 1_{\{ T^{(n)}(\Theta) \in E \}} 
\tend{\rm d}{N \to \infty} \frac{\mu(E)}{\sqrt{\pi}} \, |\cN| 
\label{eq: DKlim}
\end{align}
for all Borel set $ E $ with $ \mu(E) < \infty $. 
\end{Thm}

Theorem \ref{P DKlaw} will be proved in Section \ref{Pelikan}, 
where we adopt a finer partition for the skew-product 
than that of Theorem \ref{arclaw} 
in order to avoid a certain technical difficulty.

Second, we take up 
\begin{align}
\tau_1(x) = 
\begin{cases}
2x & (0 < x < 1/4) \\
1-2x & (1/4 < x < 1/2) \\
2-2x & (1/2 < x < 3/2) \\
2x-1 & (3/2 < x < 1) 
\end{cases}
, \quad 
\tau_2(x) = 
\begin{cases}
1-x/4 & (0 < x < 1/2) \\
(1-x)/4 & (1/2 < x < 1) 
\end{cases}
. 
\label{}
\end{align}
and we call the corresponding random map $ T $ the \emph{modified Boyarsky--G\'ora--Islam map} 
or \emph{mBGI map} in short. 
Note that both 0 and 1 are repelling fixed points for $ \tau_1 $ 
and the two-point set $ \{ 0,1 \} $ is an attractor for $ \tau_2 $. 
We choose $ p = 2/3 $ so that the two-point set $ \{ 0,1 \} $ is 
an \emph{indifferent-in-average attractor} 
in the sense that 
\begin{align}
T \{ 0,1 \} = \{ 0,1 \} , \ \bE \log |T'(0+)| = 0 , \ \bE \log |T'(1-)| = 0 . 
\label{}
\end{align}
We introduce the $ \lambda $-partition $ \alpha $ of $ [0,1] $ 
for the mBGI map as $ \alpha = \cbra{ I^-_k, I^+_k }_{k=0}^{\infty } $ with 
\begin{align}
I^-_k = \rbra{ \frac{1}{2^{k+2}}, \frac{1}{2^{k+1}} } 
, \quad 
I^+_k = \rbra{ 1-\frac{1}{2^{k+1}}, 1-\frac{1}{2^{k+2}} }. 
\label{eq: BGI alpha}
\end{align}
Although the arcsine law does not make sense, because the orbit immediately commute 
between the neighborhood of 0 and that of 1, 
we obtain the Darling--Kac law as follows.

\begin{Thm} \label{BGI DKlaw}
Regarding the mBGI map with $ p=2/3 $, 
there exists 
a unique (up to a constant multiple) $ \lambda $-a.c. $ \sigma $-finite $ T $-invariant measure, 
which is locally constant with 
\begin{align}
\mu(I^-_k) = \mu(I^+_k) 
= \frac{3 \sqrt{2}}{8} 
\rbra{ \frac{8}{3} \cdot 2^k + \frac{1}{3} \cdot (-1)^k - 1 } \cdot 2^{-k-2} 
, \quad k=0,1,2,\ldots 
\label{eq: mmm}
\end{align}
Consequently, $ \mu $ has infinite mass near 0 and 1 
and has finite mass away from 0 and 1. 
Moreover, the Darling--Kac law \eqref{eq: DKlim} holds. 
\end{Thm}

The proof of Theorem \ref{BGI DKlaw} will be proved in a similar way to Theorem \ref{arclaw}.

\subsection*{Intermittency}

The \emph{intermittency} 
has been introduced in statistical physics 
by Pomeau--Manneville \cite{MR576270} and Manneville \cite{MR592697}, 
which describes irregular transitions 
between the \emph{laminar phase} and the \emph{turbulent burst}. 
As mathematical models of intermittency, 
dynamical systems with indifferent fixed points 
have been studied by many researchers 
such as Takahashi \cite{MR948718}, 
Collet--Ferrero \cite{MR1057449}, 
Collet--Galves--Schmitt \cite{MR1185337}, 
Collet--Galves \cite{MR1239564}, 
Mori \cite{MR1247664}, 
Campanino--Isola \cite{MR1325560,MR1366535}, 
Liverani--Saussol--Vaienti \cite{MR1695915}, 
Pollicott--Yuri \cite{MR1822105}, 
Hu \cite{MR2054191} 
and Munday--Knight \cite{MR3462301}. 
The Boole and Hata maps may be regarded also as mathematical models of intermittency 
thanks to indifferent-in-average fixed point(s). 
For other studies of statistical physics about random maps, 
see Ashwin--Aston--Nicol \cite{MR1601486}, Akimoto--Aizawa \cite{MR2721930} 
and Sato--Klages \cite{MR3952110}.

\subsection*{Organization}

This paper is organized as follows. 
In Section \ref{conjugacy} we review the general theory 
assuring that a dynamical system is conjugate 
to a graph shift and to a Markov chain. 
Sections \ref{hata}, \ref{mBGI} and \ref{Pelikan} 
are devoted to the proofs 
of Theorems \ref{arclaw}, \ref{BGI DKlaw} and \ref{P DKlaw}, respectively.

\subsection*{Acknowledgements}
The authors would like to express our gratitude 
to Shafiqul Islam, Jon Aaronson and Tomoki Inoue for useful advice.
The first author thanks Yushi Nakano, Fumihiko Nakamura, Hisayoshi Toyokawa and Toru Sera 
for discussions which helped improve the earlier versions of this paper. 
He also thanks Yuzuru Sato and Takuma Akimoto 
for drawing his attention to the papers \cite{MR1601486}, \cite{MR2721930} and \cite{MR3952110}.

\section{Conjugacy induced by a Markov partition} \label{conjugacy}

\subsection{Conjugacy to a graph shift}

Let us review the general theory 
which ensures that a dynamical system with a Markov partition 
is conjugate to a graph shift. 
We complement Chapter 4 of \cite{Aa} with the proofs which were omitted there. 

Let $ (X,\lambda,T) $ be a non-singular dynamical system, i.e. 
$ (X,\cB(X),\lambda) $ is a standard measure space 
and [$ \lambda \circ T^{-1}(B) = 0 $ whenever $ \lambda(B) = 0 $]. 
Note that, for a measure $ \mu $ on $ \cB(X) $ equivalent to $ \lambda $, 
the system $ (X,\mu,T) $ is again a non-singular dynamical system. 

We call $ \alpha $ a \emph{$ \lambda $-partition} 
if $ \alpha $ is a finite or countable disjoint subfamily of $ \cB(X) $ 
satisfying 
\begin{align}
\text{[$ \textstyle \sum_{s \in \alpha } s = X $ mod $ \lambda $] 
and [$ 0 < \lambda(s) < \infty $ for all $ s \in \alpha $]}. 
\label{}
\end{align}
A $ \lambda $-partition $ \alpha $ is called a \emph{Markov partition} if 
the following conditions are satisfied: 
\begin{enumerate}

\item[(M1)] 
The $ \sigma $-field 
$ \cG := \sigma \rbra{ T^{-n}s : s \in \alpha , n \ge 0 } $ 
coincides with $ \cB(X) $ mod $ \lambda $; 

\item[(M2)] 
The map $ T : s \mapsto Ts $ is invertible $ \lambda $-a.e. 
for all $ s \in \alpha $; 

\item[(M3)] 
The forward image $ Ts $ belongs to 
$ \sigma(\alpha ) $ mod $ \lambda $ 
for all $ s \in \alpha $. 

\end{enumerate}

Regarding each element of $ \alpha $ as a point, 
we denote the space of sequences of $ \alpha $ by 
\begin{align}
Y := \cbra{ \vs = (s_0,s_1,s_2,\ldots) : 
s_0,s_1,s_2,\ldots \in \alpha } . 
\label{}
\end{align}
Let $ \theta $ denote the shift operator of $ Y $, that is, 
\begin{align}
\theta : Y \ni (s_0,s_1,s_2,\ldots) 
\mapsto (s_1,s_2,\ldots) \in Y . 
\label{}
\end{align}
A \emph{cylinder} of $ Y $ is given of the form 
\begin{align}
[s_0,s_1,\ldots,s_n] 
= \cbra{ (s_0,s_1,\ldots,s_n,t_{n+1},t_{n+2},\ldots) \in Y : 
\text{$ t_{n+1},t_{n+2},\ldots \in \alpha $} } 
\label{}
\end{align}
and we equip $ Y $ with the product topology, 
so that $ \cB(Y) $ is generated by cylinders. 
We introduce the \emph{graph shift} by 
\begin{align}
Y_\lambda := \cbra{ \vs \in Y : 
\lambda(s_0 \cap T^{-1}s_1) > 0 , \ \lambda(s_1 \cap T^{-1}s_2) > 0,\ldots } . 
\label{}
\end{align}
We note that $ \cB(Y_\lambda) $ is generated by cylinders contained in $ Y_\lambda $. 
We define a surjection $ \phi:X \to Y_\lambda $ by 
\begin{align}
\text{$ \phi(x) = \vs $ with $ x \in s_0 $, $ T(x) \in s_1 $, $ T^2(x) \in s_2 ,\ldots $} 
\label{}
\end{align}
This map $ \phi $ is Borel measurable since 
\begin{align}
\phi^{-1} [s_0,s_1,\ldots,s_n] 
= s_0 \cap T^{-1} s_1 \cap \cdots \cap T^{-n} s_n 
\in \cB(X) . 
\label{eq: phi-1}
\end{align}
We utilize the following theorem (see, e.g., \cite[Proposition 4.2.3]{Aa}) 
without proof. 

\begin{Thm} \label{thm: conjugacy}
Let $ \alpha $ be a Markov partition 
and let $ \mu $ be a measure on $ X $ which is equivalent to $ \lambda $. 
Then the Borel surjection 
$ \phi : (X,\mu,T) \to (Y_\lambda,\mu \circ \phi^{-1},\theta) $ 
is a conjugacy in the sense that 
$ \phi \circ T = \theta \circ \phi $ 
and $ \phi $ is essentially a Borel isomorphism. 
\end{Thm}

\subsection{Markov chain}

Let us study the condition 
that $ (Y_\lambda,\mu \circ \phi^{-1},\theta) $ becomes a Markov chain. 
We introduce the transition probability on $ Y $ as 
\begin{align}
p(s,t) = \frac{\lambda(s \cap T^{-1}t)}{\lambda(s)} 
, \quad s,t \in \alpha . 
\label{}
\end{align}
Note that the graph shift $ Y_\lambda $ can be rewritten as 
\begin{align}
Y_\lambda := \cbra{ \vs \in Y : 
p(s_0,s_1) > 0 , \ p(s_1,s_2) > 0,\ldots } . 
\label{}
\end{align}
Kolmogorov's extension theorem shows that, 
for any $ s \in \alpha $, 
there exists a unique probability measure $ \nu_s $ on $ Y $ 
such that 
\begin{align}
\nu_s([s_0,s_1,\ldots,s_n]) = \delta_s(s_0) 
p(s_0,s_1) p(s_1,s_2) \cdots p(s_{n-1},s_n) 
\label{}
\end{align}
for all cylinders, where $ \delta_s $ denotes the Dirac mass at $ s $. 
For a measure $ \mu $ on $ X $, 
we define a measure $ \nu_\mu $ on $ Y $ by 
\begin{align}
\nu_\mu(B) = \sum_{s \in \alpha } \mu(s) \nu_s(B) 
, \quad B \in \cB(Y) . 
\label{}
\end{align}
Note that $ \nu_{\mu} $ is supported on the graph shift $ Y_\lambda $. 

We now see, by the help of Theorem \ref{thm: conjugacy}, that the conjugacy between 
the dynamical system $ (X,\mu,T) $ and the Markov chain $ (Y_\lambda,\nu_\mu,\theta) $ 
is equivalent to the condition $ [\nu_\mu = \mu \circ \phi^{-1}] $. 
For a $ \sigma $-finite measure $ \mu $ on $ X $, 
we write $ \hat{T}_\mu : L^1(\mu) \to L^1(\mu) $ 
for the Perron--Frobenius operator of $ T $ with respect to $ \mu $: 
\begin{align}
\int_X \hat{T}_\mu f \cdot g \, \d \mu 
= \int_X f \cdot g \circ T \, \d \mu 
, \quad f \in L^1(\mu) , \ g \in L^{\infty }(\mu) . 
\label{}
\end{align}
The following proposition plays an important role in checking the condition 
$ [\nu_\mu = \mu \circ \phi^{-1}] $. 

\begin{Prop} \label{prop: numu}
Let $ \alpha $ be a Markov partition 
and let $ \mu $ be a measure on $ X $ which is equivalent to $ \lambda $. 
Then the condition $ [\nu_\mu = \mu \circ \phi^{-1}] $ is satisfied 
if and only if 
\begin{align}
\hat{T}_\mu 1_s = \sum_{t \in \alpha } c_\mu(s,t) 1_t 
\ \text{with} \ 
c_\mu(s,t) = \frac{\mu(s) p(s,t)}{\mu(t)} 
, \quad s \in \alpha 
\label{eq: hatmu}
\end{align}
holds. 
In this case, 
$ \mu $ is $ T $-invariant if and only if $ \nu_\mu $ is $ \theta $-invariant. 
\end{Prop}

\Proof{
(i) Suppose $ \nu_\mu = \mu \circ \phi^{-1} $. 
For $ s_1,\ldots,s_n \in \alpha $, we set 
\begin{align}
g = 1_{s_1 \cap T^{-1} s_2 \cap \cdots \cap T^{-(n-1)} s_n} . 
\label{}
\end{align}
We then have 
\begin{align}
& \int_X \hat{T}_\mu 1_s \cdot g \, \d \mu 
= \int_X 1_{s \cap T^{-1}s_1 \cap T^{-2} s_2 \cap \cdots \cap T^{-n} s_n} \, \d \mu 
\label{} \\
=& \mu \circ \phi^{-1}([s,s_1,\ldots,s_n]) 
= \nu_\mu ([s,s_1,\ldots,s_n]) 
= \mu(s) p(s,s_1) p(s_1,s_2) \cdots p(s_{n-1},s_n) 
\label{} \\
=& \frac{\mu(s) p(s,s_1)}{\mu(s_1)} \mu(s_1) p(s_1,s_2) \cdots p(s_{n-1},s_n) 
= \int_X \rbra{ \sum_{t \in \alpha } \frac{\mu(s) p(s,t)}{\mu(t)} 1_t } 
\cdot g \, \d \mu , 
\label{}
\end{align}
which shows \eqref{eq: hatmu}, 
since $ \cG = \cB(X) $ mod $ \mu $. 

In this case, we have 
$ \nu_\mu \circ \theta^{-1} = \mu \circ (\theta \circ \phi)^{-1} 
= \mu \circ (\phi \circ T)^{-1} = (\mu \circ T^{-1}) \circ \phi^{-1} $, 
which shows the equivalence between 
$ \mu \circ T^{-1} = \mu $ and $ \nu_\mu \circ \theta^{-1} = \nu_\mu $.

(ii) Suppose \eqref{eq: hatmu} is satisfied. 
By \eqref{eq: phi-1}, the condition $ [\nu_\mu = \mu \circ \phi^{-1}] $ 
is equivalent to 
\begin{align}
\nu_\mu([s_0,s_1,\ldots,s_n]) 
= \mu (s_0 \cap T^{-1} s_1 \cap \cdots \cap T^{-n} s_n) 
, \quad s_0,\ldots,s_n \in \alpha 
\label{eq: numu}
\end{align}
for all $ n \ge 0 $, 
since $ \cG = \cB(X) $ mod $ \mu $. 
We shall prove \eqref{eq: numu} by induction in $ n \ge 0 $. 
The case $ n=0 $ is obvious. 
Suppose \eqref{eq: numu} is true for $ n \ge 0 $. 
For $ n+1 $, we have 
\begin{align}
& \mu (s_0 \cap T^{-1} s_1 \cap \cdots \cap T^{-n} s_n \cap T^{-(n+1)} s_{n+1}) 
\label{} \\
=& \int_X 1_{s_0} \cdot 1_{s_1 \cap T^{-1} s_2 \cap \cdots \cap T^{-(n-1)} s_n \cap T^{-n} s_{n+1}} \circ T \, \d \mu 
\label{} \\
=& \int_X \hat{T}_\mu 1_{s_0} \cdot 1_{s_1 \cap T^{-1} s_2 \cap \cdots \cap T^{-(n-1)} s_n \cap T^{-n} s_{n+1}} \, \d \mu 
\label{} \\
=& \frac{\mu(s_0) p(s_0,s_1)}{\mu(s_1)} 
\mu(s_1 \cap T^{-1}s_2 \cap \cdots \cap T^{-(n-1)} s_n \cap T^{-n} s_{n+1}) . 
\label{eq: conj eq}
\end{align}
By the induction assumption, we have 
\begin{align}
\eqref{eq: conj eq} 
=& \frac{\mu(s_0) p(s_0,s_1)}{\mu(s_1)} \nu_\mu([s_1,s_2,\ldots,s_n,s_{n+1}]) 
\label{} \\
=& \mu(s_0) p(s_0,s_1) p(s_1,s_2) \cdots p(s_n,s_{n+1}) 
\label{} \\
=& \nu_\mu([s_0,s_1,\ldots,s_n,s_{n+1}]) , 
\label{}
\end{align}
which shows \eqref{eq: numu} for $ n+1 $. 
}

\subsection{Conservative ergodicity} \label{sec: CE}

It is well-known 
(see \cite[Proposition 1.2.2]{Aa}) 
that the dynamical system $ (X,\lambda,T) $ is \emph{conservative ergodic} if and only if 
\begin{align}
\sum_{n=1}^{\infty } 1_B \circ T^n = \infty 
\quad \text{$ \lambda $-a.e. for all $ B \in \cB(X) $ with $ \lambda(B)>0 $} . 
\label{}
\end{align}
This property can be discussed in terms of the Markov chain. 

We say that a transition matrix $ P = (p(s,t))_{s,t \in \alpha } $ is \emph{irreducible} if 
for any $ s,t \in \alpha $ there exists $ n \ge 1 $ such that 
\begin{align}
p^{(n)}(s,t) := \sum_{s_1,\ldots,s_{n-1} \in \alpha } p(s,s_1) p(s_1,s_2) \cdots p(s_{n-2},s_{n-1}) p(s_{n-1},t) > 0 
\label{}
\end{align}
(we understand that $ p^{(1)}(s,t) = p(s,t) $). 
We say that an irreducible transition matrix $ P $ 
is \emph{recurrent} (resp. \emph{transient}) if 
\begin{align}
\sum_{n=1}^{\infty } p^{(n)}(s,s) = \infty \ (\text{resp. $ < \infty $}) 
\quad \text{for some $ s \in \alpha $}, 
\label{eq: ir}
\end{align}
and in this case \eqref{eq: ir} holds for all $ s \in \alpha $. 
If we denote the first hitting time of $ t \in \alpha $ by 
\begin{align}
\varphi_t(\vs) = \inf \{ n \ge 1 : s_n = t \} 
\quad (\vs = (s_0,s_1,s_2,\ldots) \in Y), 
\label{}
\end{align}
then it is well-known (see, e.g. \cite[Theorems 1.5.3 and 1.5.7]{No}) that 
the condition \eqref{eq: ir} can be replaced by 
\begin{align}
\nu_s(\varphi_s < \infty ) = 1 \ (\text{resp. $ < 1 $}) 
\quad \text{for some $ s \in \alpha $} . 
\label{}
\end{align}
It is also well-known (see e.g. \cite[Theorems 1.7.5 and 1.7.6]{No}) that, 
if the transition matrix $ P $ is irreducible recurrent, 
then 
\begin{align}
\nu_s(\varphi_t < \infty ) = 1 
\quad \text{for all $ s,t \in \alpha $} , 
\label{}
\end{align}
and there exists a unique (up to a constant multiple) 
$ \sigma $-finite $ P $-invariant measure on $ \alpha $: 
\begin{align}
\rho(t) = \sum_{s \in \alpha } \rho(s) p(s,t) 
, \quad t \in \alpha . 
\label{}
\end{align}

\begin{Prop} \label{prop: P-inv ce} 
Suppose $ [\nu_\mu = \mu \circ \phi^{-1}] $ is satisfied. 
Suppose, in addition, the transition matrix $ P $ is irreducible recurrent. 
Let $ \rho $ be a unique (up to a constant multiple) 
$ \sigma $-finite $ P $-invariant measure on $ \alpha $. 
Define the measure $ \mu $ on $ X $ by 
\begin{align}
\mu(\d x) = \sum_{s \in \alpha } \frac{\rho(s)}{\lambda(s)} 1_s(x) \d x . 
\label{}
\end{align}
Then it holds that the dynamical system $ (X,\mu,T) $ is conservative ergodic, 
and that the measure $ \mu $ is 
a unique (up to a constant multiple) $ \lambda $-a.c. $ \sigma $-finite $ T $-invariant measure. 
\end{Prop}

\Proof{
By \cite[Theorem 4.5.3]{Aa}, we see that 
the dynamical system $ (Y_\lambda,\nu_\mu,\theta) $ is conservative ergodic. 
Hence, so is $ (X,\mu,T) $, by the conjugacy. 

The $ T $-invariance of $ \mu $ follows from the $ \theta $-invariance of $ \nu_\mu $, 
by the conjugacy. 
The uniqueness follows from \cite[Theorem 1.5.6]{Aa}. 
}

\section{Proof for the Hata map} \label{hata}

For the Hata map \eqref{eq: hatamap}, 
we define the transition probability on $ \alpha $ as 
\begin{align}
q(s,t) 
= \frac{\bE \lambda(s \cap T^{-1} t)}{\lambda(s)} 
= \frac{p \lambda(s \cap \tau_1^{-1} t) + (1-p) \lambda(s \cap \tau_2^{-1} t)}
{\lambda(s)} 
, \quad s,t \in \alpha . 
\label{eq: prmt}
\end{align}

\begin{Prop} \label{prop: Hata P} 
Regarding the Hata map with $ p=1/2 $, 
the transition matrix $ Q=(q(s,t))_{s,t \in \alpha } $ is irreducible recurrent 
and has a unique (up to a constant multiple) $ Q $-invariant measure $ \rho $ given as 
\begin{align}
\rho(I^-_k) = \rho(I^+_k) 
= (2^{k+1}-1) \cdot 2^{-k-2} 
, \quad k=0,1,2,\ldots , 
\label{eq: rhoim2}
\end{align}
where $ I^-_k $'s and $ I^+_k $'s have been introduced in \eqref{eq: alpha}. 
\end{Prop}

\Proof{
We sometimes omit ``mod $ \lambda $" 
in the identities among subsets of $ X $. 
Note that 
\begin{align}
\tau_1 I^-_k = I^-_{k+1} \ (k \ge 0) 
, \quad 
\tau_2 I^-_k = 
\begin{cases}
I^-_{k-1} & (k \ge 1) 
\\
\sum_{j=0}^{\infty } I^+_j & (k=0) 
\end{cases}
, 
\label{} \\
\tau_1 I^+_k = 
\begin{cases}
I^+_{k-1} & (k \ge 1) 
\\
\sum_{j=0}^{\infty } I^-_j & (k=0) 
\end{cases}
, \quad 
\tau_2 I^+_k = I^+_{k+1} \ (k \ge 0) 
. 
\label{}
\end{align}
After an easy computation we have 
\begin{align}
q(s,t) = 
\begin{cases}
\frac{1}{2} & \text{if $ (s,t) = (I^-_k,I^-_{k+1}) $ or $ (I^+_k,I^+_{k+1}) $ for $ k \ge 0 $} 
\\
\frac{1}{2} & \text{if $ (s,t) = (I^-_k,I^-_{k-1}) $ or $ (I^+_k,I^+_{k-1}) $ for $ k \ge 1 $} 
\\
\frac{1}{2^{j+2}} & \text{if $ (s,t) = (I^-_0,I^+_j) $ or $ (I^+_0,I^-_j) $ for $ j \ge 0 $} 
\end{cases}
\label{}
\end{align}
(Note that $ \frac{1}{2} + \sum_{j=0}^{\infty } \frac{1}{2^{j+2}} = 1 $). 
It is obvious that $ Q $ is irreducible. 
So we need only to show $ \nu_{I^+_0}(\varphi_{I^+_0}<\infty )=1 $ 
(see Section \ref{sec: CE}). 

Since the simple symmetric random walk on $ \bZ $ is irreducible recurrent, we have 
\begin{align}
\nu_{I^-_k}(\varphi_{I^-_0}<\infty ) = 1 
, \quad k \ge 1 . 
\label{}
\end{align}
By the Markov property, we have 
\begin{align}
\nu_{I^+_0}(\varphi_{I^-_0}<\infty ) 
= q(I^+_0,I^-_0) + \sum_{k \ge 1} q(I^+_0,I^-_k) \nu_{I^-_k}(\varphi_{I^-_0} < \infty ) 
= 1 . 
\label{}
\end{align}
We also have $ \nu_{I^-_0}(\varphi_{I^+_0}<\infty )=1 $ by symmetry. 
By the strong Markov property, we obtain 
\begin{align}
\nu_{I^+_0}(\varphi_{I^+_0}<\infty ) 
\ge \nu_{I^+_0}(\varphi_{I^-_0}<\infty ) \nu_{I^-_0}(\varphi_{I^+_0}<\infty ) = 1 , 
\label{}
\end{align}
which shows $ \nu_{I^+_0}(\varphi_{I^+_0}<\infty ) = 1 $, 
and hence $ Q $ is recurrent. 

Let us prove that 
the measure $ \rho $ given in \eqref{eq: rhoim2} is $ Q $-invariant. 
By symmetry, the $ Q $-invariance is equivalent to the recurrence relation: 
\begin{align}
\rho(I^-_k) 
= \frac{1}{2} \rho(I^+_{k+1}) + \frac{1}{2} \rho(I^+_{k-1}) + \frac{1}{2^{k+2}} \rho(I^-_0) 
, \quad k=0,1,2,\ldots , 
\label{}
\end{align}
if we understand that $ \rho(I^-_{-1}) = 0 $. 
We can easily verify that the measure $ \rho $ given in \eqref{eq: rhoim2} satisfies 
this recurrence relation. 
}

Let us represent the random map $ T $ of \eqref{eq: randommap} by the skew-product. 
In what follows, we write $ \lambda $ for the Lebesgue measure on $ X = [0,1] $ 
and write $ (\Omega, \cB(\Omega), \bP) $ 
for the coin flip: 
\begin{align}
\Omega = \cbra{ \omega = (\omega_1,\omega_2,\ldots) : 
\omega_1,\omega_2,\ldots \in \{ \tau_1,\tau_2 \} } 
\label{}
\end{align}
equipped with the product topology, 
$ \cB(\Omega) = $ the Borel field of $ \Omega $, and 
\begin{align}
\bP(\omega \in \Omega : 
\omega_1=\chi_1,\ldots,\omega_n=\chi_n) = p^{N^1_n(\chi)} (1-p)^{N^2_n(\chi)} 
\label{}
\end{align}
for all $ \chi_1,\ldots,\chi_n \in \{ \tau_1,\tau_2 \} $ and all $ n $, 
where $ N^1_n(\chi) = \# \{ k=1,\ldots,n : \chi_k = \tau_1 \} $ 
and $ N^2_n(\chi) = \# \{ k=1,\ldots,n : \chi_k = \tau_2 \} $. 
For $ B \in \cB(X) $, we denote $ \tilde{B} = B^{\sim} = \Omega \times B $, 
and we equip $ \tilde{X} = \Omega \times X $ 
with the product topology of $ \Omega $ and $ X $. 
We write $ \cB(\tilde{X}) $ for the Borel field of $ \tilde{X} $ 
and write $ \tilde{\lambda} = \bP \otimes \lambda $ 
for the product measure of $ \bP $ and $ \lambda $. 
We define the deterministic transformation $ \tilde{T} : \tilde{X} \to \tilde{X} $ as 
\begin{equation}
\tilde{T}(\omega,x) = \rbra{ \theta \omega, \omega_1 x } 
, \quad 
\omega = (\omega_1, \omega_2, \ldots) \in \Omega, \ x \in X . 
\end{equation}
Here $ \theta : \Omega \to \Omega $ stands for the shift operator: 
$ \theta(\omega_1,\omega_2,\ldots) = (\omega_2,\omega_3,\ldots) $.
The random maps $ \{ T_n \} $ and $ \{ T^{(n)} \} $ in Section \ref{sec: intro} 
are obtained as 
\begin{align}
T_n(\omega)(x) = \omega_n(x) 
, \quad 
T^{(n)}(\omega)(x) = \omega_n \omega_{n-1} \cdots \omega_1(x) , 
\label{}
\end{align}
so that $ \tilde{T}^n(\omega,x) = (\theta_n \omega,T^{(n)}(\omega)(x)) $. 
For a measure $ \mu $ on $ X $, 
it is easy to see that $ \tilde{\mu} := \bP \otimes \mu $ is $ \tilde{T} $-invariant 
if and only if $ \mu $ is $ T $-invariant in the sense of \eqref{eq: T-inv}.

\begin{Prop} \label{Markovpart}
Regarding the Hata map with $ p=1/2 $, 
the family 
\begin{align}
\tilde{\alpha} 
= \cbra{ \tilde{s} : s \in \alpha } 
= \cbra{ \tilde{I}^-_k, \tilde{I}^+_k }_{k=0}^{\infty } 
\label{}
\end{align}
is a Markov partition for the dynamical system 
$ (\tilde{X},\tilde{\mu},\tilde{T}) $, 
where $ \mu $ is given as \eqref{eq: im2}. 
\end{Prop}

\Proof{
By Proposition \ref{prop: Hata P}, we see that 
$ \tilde{\mu} = \bP \otimes \mu $ is $ \tilde{T} $-invariant, 
so that $ \mu $ is $ T $-invariant. 

We sometimes omit ``mod $ \tilde{\mu} $" 
in the identities among subsets of $ \tilde{X} $. 
Let us prove (M2) and (M3) at once. 
For $ (\omega,x) \in \tilde{I}^-_k \in \tilde{\alpha } $ with $ k \ge 1 $, we have 
\begin{align}
\text{$ \tilde{T}(\omega,x) = (\theta \omega,\omega_1(x)) 
\in \tilde{I}^-_{k+1} $ or $ \tilde{I}^-_{k-1} $ 
according as $ \omega_1 = \tau_1 $ or $ \tau_2 $}. 
\label{}
\end{align}
We thus see that $ \tilde{T} : \tilde{I}^-_k \to \tilde{T}(\tilde{I}^-_k) 
= \tilde{I}^-_{k+1} + \tilde{I}^-_{k-1} \in \sigma( \tilde{\alpha } ) $ 
is invertible as 
\begin{align}
\tilde{T}^{-1}(\omega,x) = 
\begin{cases}
((\tau_1,\omega),2x) & ((\omega,x) \in \tilde{I}^-_{k+1}) , 
\\
((\tau_2,\omega),x/2) & ((\omega,x) \in \tilde{I}^-_{k-1}) . 
\end{cases}
\label{}
\end{align}
For $ (\omega,x) \in \tilde{I}^-_0 \in \tilde{\alpha } $, we have 
\begin{align}
\text{$ \tilde{T}(\omega,x) = (\theta \omega,\omega_1(x)) 
\in \tilde{I}^-_1 $ or $ (1/2,1)^{\sim} $ according as $ \omega_1 = \tau_1 $ or $ \tau_2 $}. 
\label{}
\end{align}
Since $ \tilde{T}(\tilde{I}^-_0) 
= \tilde{I}^-_1 + (1/2,1)^{\sim} 
= \tilde{I}^-_1 + \sum_{k=0}^{\infty } \tilde{I}^+_k \in \sigma( \tilde{\alpha } ) $, 
we thus see that $ \tilde{T} : \tilde{I}^-_0 \to \tilde{T}(\tilde{I}^-_0) $ 
is invertible as 
\begin{align}
\tilde{T}^{-1}(\omega,x) = 
\begin{cases}
((\tau_1,\omega),2x) & ((\omega ,x) \in \tilde{I}^-_1) , 
\\
((\tau_2,\omega),x/2) & ((\omega ,x) \in \tilde{I}^+_k, \ k \ge 0) . 
\end{cases}
\label{}
\end{align}
By symmetry we obtain a similar result 
for $ \tilde{T} : \tilde{I}^+_k \to \tilde{T}(\tilde{I}^+_k) $. 

Let us now prove (M1). 
For $ \omega = (\omega_1,\omega_2,\ldots) $ and for $ x \in X $, 
we write $ \xi_n(\omega,x) = \omega_n $. 
It suffices to show 
\begin{align}
& \cbra{ \xi_n = \tau_1 }, \ \cbra{ \xi_n = \tau_2 } \in \cG 
, \quad n=1,2,\ldots, 
\label{eq: step1} \\
& \rbra{ \frac{i}{2^m}, \frac{i+1}{2^m} }^{\sim} \in \cG 
, \quad i=0,1,\ldots,2^m-1, \ m=0,1,2,\ldots . 
\label{eq: step2}
\end{align}
To obtain $ \{ \xi_1 = \tau_1 \} \in \cG $, we have 
\begin{align}
\cbra{ \xi_1 = \tau_1 } 
=& \sum_{k \geq 0} \cbra{ \xi_1 = \tau_1 } \cap \tilde{I}^-_k 
+ \sum_{k \geq 0} \cbra{ \xi_1 = \tau_1 } \cap \tilde{I}^+_k 
\label{} \\
=& \sum_{k \geq 0} \rbra{ \tilde{T}^{-1} \tilde{I}^-_{k+1} } \cap \tilde{I}^-_{k} + \sum_{k \geq 1} \rbra{ \tilde{T}^{-1} \tilde{I}^+_{k-1} } \cap \tilde{I}^+_{k} + \rbra{ \sum_{j \geq 0} \tilde{T}^{-1} \tilde{I}^-_{j} } \cap \tilde{I}^+_{0} \in \cG, 
\label{}
\end{align}
and hence we obtain 
\begin{align}
\cbra{ \xi_n = \tau_1 } = \tilde{T}^{-(n-1)} \cbra{ \xi_1 = \tau_1 } \in \mathcal{G}. 
\label{}
\end{align}
Since $ \cbra{ \xi_n = \tau_2 } = \cbra{ \xi_n = \tau_1 }^c \in \cG $, 
we have proved \eqref{eq: step1}. 

We proceed to show \eqref{eq: step2} by induction. 
The case $ m=0 $ is obvious; $ (0,1)^{\sim} = \tilde{X} \in \cG $. 
Suppose \eqref{eq: step2} holds for $ m $ 
and let us prove that \eqref{eq: step2} holds also for $ m+1 $. 
Set 
\begin{align}
\eta_n = 
\begin{cases}
1 & (\xi_n=\tau_1) \\
-1 & (\xi_n=\tau_2) 
\end{cases}
\quad (n \ge 1) 
\label{}
\end{align}
and 
\begin{align}
W_n = \eta_1 + \cdots + \eta_n , \quad n \ge 1 . 
\label{eq: W_n}
\end{align}
Note that $ (\{ W_n \},\bP) $ is the simple symmetric random walk on $ \bZ $ starting from 0. 
For $ k \in \bZ $, we denote the first hitting time of $ k $ for the random walk by 
\begin{align}
\varphi^W_k = \inf \{ n \ge 1 : W_n = k  \} . 
\label{eq: varphi^W_a}
\end{align}
Since $ (\{ W_n \},\bP) $ is irreducible recurrent, we have 
$ \bP(\varphi^W_k < \infty ) = 1 $ for all $ k \in \bZ $. 

For any $ i=0,\ldots,2^m-1 $ so that $ \frac{i+1}{2^{m+1}} \leq \frac{1}{2} $, 
we want to show $ \Big. \rbra{ \frac{i}{2^{m+1}}, \frac{i+1}{2^{m+1}} }^{\sim} \in \cG $. 
We divide it into 
\begin{align}
\rbra{ \frac{i}{2^{m+1}}, \frac{i+1}{2^{m+1}} }^{\sim} 
=& \sum_{n=1}^{\infty } \rbra{ \frac{i}{2^{m+1}}, \frac{i+1}{2^{m+1}} }^{\sim} 
\cap \cbra{ \varphi^W_{-1} =n } , 
\label{}
\end{align}
where $ \varphi^W_{-1} $ can be regarded as the first hitting time of $ m $ 
for the random walk starting from $ m+1 $. 
On the event $ \{ \varphi^W_{-1} = n \} $, 
we have $ \xi_{j-1} \circ \cdots \circ \xi_1 (x) \in \rbra{ 0,\frac{i+1}{2^{m+1}} } $ 
for all $ j \le n $, so that $ \xi_j \circ \xi_{j-1} \circ \cdots \circ \xi_1 (x) 
= 2^{-W_j} x $ for all $ j \le n $. 
Hence 
\begin{align}
\rbra{ \frac{i}{2^{m+1}}, \frac{i+1}{2^{m+1}} }^{\sim} \cap \cbra{ \varphi^W_{-1} =n } 
= \rbra{ \tilde{T}^{-n} \rbra{ \frac{i}{2^{m}}, \frac{i+1}{2^{m}}}^{\sim} } 
\cap \cbra{ \varphi^W_{-1}=n } . 
\label{}
\end{align}
Since $ \tilde{T}^{-n} \rbra{ \frac{i}{2^{m}}, \frac{i+1}{2^{m}}}^{\sim} \in \cG $ 
by the assumption of induction 
and since 
\begin{align}
\cbra{ \varphi^W_{-1}=n } \in \sigma(\xi_1,\ldots,\xi_n) \in \cG , 
\label{}
\end{align}
we obtain $ \Big. \rbra{ \frac{i}{2^{m+1}}, \frac{i+1}{2^{m+1}} }^{\sim} \in \cG $ 
for $ i=0,1,\ldots,2^m-1 $. 
By symmetry, we also obtain 
$ \Big. \rbra{ \frac{i}{2^{m+1}}, \frac{i+1}{2^{m+1}} }^{\sim} \in \cG $ 
for $ i=2^m,2^m+1,\ldots,2^{m+1}-1 $. 
We have now obtained \eqref{eq: step2}. 

The proof is therefore complete.
}

We have seen in Proposition \ref{Markovpart} 
that $ \tilde{\alpha } $ is a Markov partition, 
and so we see by Theorem \ref{thm: conjugacy} 
that our dynamical system $ (\tilde{X},\tilde{\mu},\tilde{T}) $ 
is conjugate to $ (\tilde{Y}_{\tilde{\lambda}},\tilde{\mu} \circ \tilde{\phi}^{-1},\theta) $, 
where $ \tilde{Y} $, $ \tilde{\phi} $, etc. 
are so defined as in Section \ref{conjugacy}. 
Note that, for $ s,t \in \alpha $, we have 
\begin{align}
\tilde{p}(\tilde{s},\tilde{t}) 
= \frac{\tilde{\lambda} \bigl( \tilde{s} \cap \tilde{T}^{-1} \tilde{t} \bigr)}{\tilde{\lambda}(\tilde{s})} 
= \frac{\bE \lambda(s \cap T^{-1} t)}{\lambda(s)} 
= q(s,t) , 
\label{}
\end{align}
and so the transition matrix $ \tilde{P} = (\tilde{p}(\tilde{s},\tilde{t}))_{\tilde{s},\tilde{t} \in \tilde{\alpha }} $ is irreducible recurrent.

\begin{Prop} \label{thm: Hata ce}
Regarding the Hata map with $ p=1/2 $, 
the dynamical system $ (\tilde{X},\tilde{\mu},\tilde{T}) $ 
is conjugate to the Markov chain $ (\tilde{X},\tilde{\nu}_{\tilde{\mu}},\tilde{T}) $, 
which is irreducible recurrent. 
Consequently, the dynamical system $ (\tilde{X},\tilde{\lambda},\tilde{T}) $ 
is conservative ergodic, and so the measure $ \tilde{\mu} $ 
is a unique (up to a constant multiple) 
$ \tilde{\lambda} $-a.c. $ \sigma $-finite $ \tilde{T} $-invariant measure. 
\end{Prop}

\Proof{
We write $ \tilde{T}^{\wedge}_{\tilde{\mu}} $ 
for the Perron--Frobenius operator of $ \tilde{T} $ with respect to $ \tilde{\mu} $. 
For $ A \in \cB(\Omega) $, $ B \in \cB(X) $ and $ s \in \alpha $, we have 
\begin{align}
& \int_{\tilde{X}} 
\tilde{T}^{\wedge}_{\tilde{\mu}} 1_{\tilde{s}} \cdot 1_{A \times B} \, \d \tilde{\mu} 
= \int_{\tilde{X}} 1_{\tilde{s}} \cdot 1_{A \times B} \circ \tilde{T} \, \d \tilde{\mu} 
\label{} \\
=& \frac{\mu(s)}{\lambda(s)} 
\int_{\Omega} \int_X 1_s(x) \cdot 1_A(\theta \omega) 1_B(\omega_1(x)) \, \bP(\d \omega) \d x 
\label{} \\
=& \frac{\mu(s)}{\lambda(s)} \bP(A) \cdot \frac{1}{2} \sum_{i=1,2} \int_X 
1_s(\tau_i^{-1}(x)) \cdot 1_B(x) \, (\tau_i^{-1})' \d x 
\label{} \\
=& \frac{\mu(s)}{\lambda(s)} \bP(A) \cdot \frac{1}{2} \sum_{i=1,2} \int_X \sum_{t \in \alpha } 
1_{\tau_i s \cap t}(x) \cdot 1_B(x) \, (\tau_i^{-1})' \d x 
\label{} \\
=& \int_{\tilde{X}} \sum_{t \in \alpha } \rbra{ \frac{\mu(s)}{\lambda(s)} \cdot \frac{\lambda(t)}{\mu(t)} \cdot \frac{1}{2} \sum_{i=1,2} 
(\tau_i^{-1})' 1_{(\tau_i s \cap t)^{\sim}} } \cdot 1_{A \times B} \, \d \tilde{\mu} . 
\label{}
\end{align}
Since $ (\tau_i s \cap t)^{\sim} = \tilde{t} $ or $ \emptyset $ 
(see (iii) of Proposition \ref{Markovpart}) 
and since $ (\tau_i^{-1})' $ is constant on $ t $, 
we obtain the representation 
\begin{align}
\tilde{T}^{\wedge}_{\tilde{\lambda}} 1_{\tilde{s}} 
= \sum_{t \in \alpha } c(s,t) 1_{\tilde{t}} 
, \quad s \in \alpha 
\label{}
\end{align}
for some function $ c(s,t) $ on $ \alpha \times \alpha $. 
Since $ \tilde{\mu} $ is locally constant on $ \tilde{\alpha } $, we have 
\begin{align}
c(s,t) 
= \frac{1}{\tilde{\mu}(\tilde{t})} \int_{\tilde{X}}  
\tilde{T}^{\wedge}_{\tilde{\mu}} 1_{\tilde{s}} \cdot 1_{\tilde{t}} \, \d \tilde{\mu} 
= \frac{\tilde{\mu} \bigl( \tilde{s} \cap \tilde{T}^{-1} \tilde{t} \bigr)}{\tilde{\mu}(\tilde{t})}
= \frac{\tilde{\mu}(s)}{\tilde{\mu}(\tilde{t})} \cdot \frac{\tilde{\lambda} \bigl( \tilde{s} \cap \tilde{T}^{-1} \tilde{t} \bigr)}{\tilde{\lambda}(\tilde{s})} 
= \frac{\tilde{\mu}(s) \tilde{p}\bigl( \tilde{s} \cap \tilde{T}^{-1} \tilde{t} \bigr)}{\tilde{\mu}(\tilde{t})} , 
\label{}
\end{align}
which shows, by Proposition \ref{prop: numu}, 
that $ \tilde{\nu}_{\tilde{\mu}} = \tilde{\mu} \circ \tilde{\phi}^{-1} $. 
}

For the proof of Theorem \ref{arclaw}, 
we appeal to Thaler--Zweim\"uller's result \cite{TZ}. 
We set 
\begin{align}
J := \rbra{ \frac{1}{4},\frac{3}{4} } = I^-_0 + I^+_0 
, \quad 
R^- := \rbra{ 0,\frac{1}{4} } = \sum_{k=1}^{\infty } I^-_k 
, \quad 
R^+ := \rbra{ \frac{3}{4},1 } = \sum_{k=1}^{\infty } I^+_k 
\label{}
\end{align}
so that $ J + R^- + R^+ = X $, where the identities hold mod $ \lambda $. 
By the fact that 
\begin{align}
\tau_1^{-1} R^- = \rbra{ 0,\frac{1}{8} } \subset R^- 
, \quad 
\tau_2^{-1} R^- = \rbra{ 0,\frac{5}{8} } \subset J+R^- 
\quad \text{mod $ \lambda $} 
\label{}
\end{align}
and by symmetry, we see that $ \tilde{J} $ dynamically separates 
$ \tilde{R}^- $ and $ \tilde{R}^+ $, 
in the sense that 
\begin{align}
\text{$ (\omega,x) \in \tilde{R}^{\pm} $ and $ \tilde{T}^2(\omega,x) \in \tilde{R}^{\mp} $ 
imply $ \tilde{T}(\omega,x) \in \tilde{J} $}. 
\label{}
\end{align}
We may call $ J $ the \emph{junction} and $ R^+ $ and $ R^- $ the \emph{rays}. 
We denote the first return time by 
\begin{align}
\varphi_{\tilde{J}}(\omega,x) 
:= \inf \cbra{ n \ge 1 : \tilde{T}^n(\omega,x) \in \tilde{J} } 
\label{eq: tilde J}
\end{align}
and set 
\begin{align}
\tilde{J}_0 := \tilde{J} 
, \quad 
\tilde{J}_n := \tilde{J}^c \cap \cbra{ \varphi_{\tilde{J}} =n } 
\ (n=1,2,\ldots) . 
\label{eq: tilde Jn}
\end{align}
(Note that $ \tilde{J}_n $ is \emph{not} of the form $ \tilde{J}_n = \Omega \times J_n $.) 
For $ N=1,2,\ldots $, we denote the wandering rate of $ \tilde{J} $ by 
\begin{align}
w_N(\tilde{J}) 
:= \sum_{n=0}^{N-1} \tilde{\mu}(\tilde{J} \cap \{ \varphi_{\tilde{J}} > n \}) 
= \int_{\tilde{J}} \rbra{ \sum_{n=0}^{N-1} \tilde{T}^{\wedge n}_{\tilde{\mu}} 1_{\tilde{J}_n} } \d \tilde{\mu} 
\label{eq: wN tilde J}
\end{align}
and the wandering rate of $ \tilde{J} $ through $ \tilde{R}^{\pm} $ by 
\begin{align}
w_N(\tilde{J},\tilde{R}^{\pm}) 
:= \sum_{n=0}^{N-1} \tilde{\mu}(\tilde{J} \cap \tilde{T}^{-1} \tilde{R}^{\pm} 
\cap \{ \varphi_{\tilde{J}} > n \}) 
= \tilde{\mu}(\tilde{J} \cap \tilde{T}^{-1} \tilde{R}^{\pm}) 
+ \sum_{n=1}^{N-1} \tilde{\mu}(\tilde{J}_n \cap \tilde{R}^{\pm}) . 
\label{}
\end{align}

Recall that $ \varphi^W_k $ is the first hitting time of $ k $ 
for the simple symmetric random walk $ (\{ W_n \},\bP) $, 
which has been introduced in \eqref{eq: varphi^W_a}. 

\begin{Lem} \label{nthPF}
For $ n \ge 1 $, it holds that 
\begin{align}
\tilde{T}^{\wedge n}_{\tilde{\mu}} 1_{\tilde{J}_n \cap \tilde{R}^{\pm}}(\omega,x) 
= c_n 1_{I^{\pm}_0}(x) 
, \quad 
\tilde{T}^{\wedge n}_{\tilde{\mu}} 1_{\tilde{J}_n}(\omega,x) = c_n 1_J(x) 
, \quad \text{$ \tilde{\lambda} $-a.e.}, 
\label{eq: nthPF} 
\end{align}
where 
\begin{align}
c_n := \sum_{k=1}^{\infty } \rbra{ 2^{k+1}-1 } 2^{-k} \bP(\varphi^W_{-k}=n) . 
\label{eq: cn}
\end{align}
\end{Lem}

\Proof{
By symmetry and by $ \tilde{J}_n \subset \tilde{R}^+ + \tilde{R}^- $, 
it suffices to show 
\begin{align}
\int_{\tilde{X}} g_1(\omega) g_2(x) \tilde{T}^{\wedge n}_{\tilde{\mu}} 
1_{\tilde{J}_n \cap \tilde{R}^-}(\omega, x) \, \bP(\d \omega) \mu(\d x) 
= c_n \bE[g_1] \int_{I^-_0} g_2(y) \mu(\d y) 
\label{eq: T wedge n JnR-}
\end{align}
for all $ g_1 \in L^{\infty }(\bP) $ and $ g_2 \in L^{\infty }(\mu) $. 
For $ k \ge 1 $, we have 
\begin{align}
\tilde{J}_n \cap \tilde{I}^-_k 
= \tilde{I}^-_k \cap \tilde{T}^{-1} \tilde{R}^- 
\cap \cdots \cap \tilde{T}^{-(n-1)} \tilde{R}^- 
\cap \tilde{T}^{-n} \tilde{J} 
= \tilde{I}^-_k \cap \{ \varphi^W_{-k} = n \} , 
\label{}
\end{align}
where $ \varphi^W_{-k} $ can be regarded as the first hitting time of $ 0 $ 
for the random walk starting from $ k $. 
The left hand side of \eqref{eq: T wedge n JnR-} equals to 
\begin{align}
& \int_{\tilde{X}} g_1(\theta^n \omega) 
g_2(\omega_n \circ \omega_{n-1} \circ \cdots \circ \omega_1 (x)) 
1_{\tilde{J}_n}(\omega, x)1_{R^-}(x) \, \bP(\d \omega) \mu(\d x) 
\label{} \\
=& \bE[g_1] \sum_{k=1}^{\infty } \int_{\tilde{X}} 
g_2(\omega_n \circ \omega_{n-1} \circ \cdots \circ \omega_1 (x)) 
1_{\{ \varphi^W_{-k} = n \}}(\omega) 1_{I^-_k}(x) \, \bP(\d \omega) \mu(\d x) 
\label{} \\
=& \bE[g_1] \sum_{k=1}^{\infty } \int_X g_2(2^k x) 
1_{I^-_k}(x) \, \mu(\d x) \, \bP(\varphi^W_{-k}=n) 
\label{eq: f22kx} \\
=& \bE[g_1] \sum_{k=1}^{\infty } \int_X g_2(y) 
1_{I^-_0}(y) \frac{\mu(I^-_k)}{\lambda(I^-_k)} \cdot \frac{\lambda(I^-_0)}{\mu(I^-_0)} 
\cdot 2^{-k} \, \mu(\d y) \, \bP(\varphi^W_{-k}=n) , 
\label{}
\end{align}
which equals to the right hand side of \eqref{eq: T wedge n JnR-}. 
Here, the equality \eqref{eq: f22kx} can be obtained 
by the same argument as that of the proof of \eqref{eq: step2} of Proposition \ref{Markovpart}. 
}

We study the wandering rate asymptotics as follows. 

\begin{Lem} \label{lem: wandering} 
The following assertions hold as $ N \to \infty $: 
\begin{enumerate}

\item 
$ \displaystyle w_N(\tilde{J}) \sim \sqrt{\frac{2}{\pi}}N^{1/2} $ 
and 
$ \displaystyle w_N(\tilde{J}, \tilde{R}^{\pm}) \sim \frac{1}{2} \sqrt{\frac{2}{\pi}}N^{1/2} $; 

\item 
$ \displaystyle \frac{1}{w_N(\tilde{J})} 
\sum_{n=0}^{N-1} \tilde{T}^{\wedge n} 1_{\tilde{J}_n}(\omega, x) 
\to 2 \cdot 1_J(x) $ 
uniformly on $ \tilde{J} $; 

\item 
$ \displaystyle \frac{1}{w_N(\tilde{J}, \tilde{R}^{\pm})} \sum_{n=0}^{N-1} \tilde{T}^{\wedge n} 1_{\tilde{J}_n \cap \tilde{R}^{\pm}}(\omega, x) \to 4 \cdot 1_{I^{\pm}_0}(x) $ 
uniformly on $ \tilde{J} $. 

\end{enumerate}
Here by $ a_N \sim b_N $ we mean $ a_N/b_N \to 1 $. 
\end{Lem}

\Proof{
Let $ 0 < z < 1 $. By the strong Markov property and stationarity, we have 
\begin{align}
\bE z^{\varphi^W_{-1}} 
= \frac{z}{2} \bE z^{\varphi^W_{-2}} + \frac{z}{2} 
= \frac{z}{2} (\bE z^{\varphi^W_{-1}})^2 + \frac{z}{2} , 
\label{}
\end{align}
and so we obtain, for any $ k \ge 1 $, 
\begin{align}
\bE z^{\varphi^W_{-k}} 
= \rbra{ \bE z^{\varphi^W_{-1}} }^{k} 
= \rbra{ \frac{1 - \sqrt{1-z^2}}{z} }^{k} . 
\label{}
\end{align}

Let us study the asymptotic behavior. As $ z \up 1 $, we have 
\begin{align}
w := 1 - \bE z^{\varphi^W_{-1}} 
= \frac{\sqrt{1-z}(\sqrt{1+z} - \sqrt{1-z})}{z} 
\sim \sqrt{2} \sqrt{1-z} . 
\label{eq: Hata w}
\end{align}
The generating function of $ \{ c_n \} $ can be computed as 
\begin{align}
\sum_{n=1}^{\infty } c_n z^n 
=& \sum_{k=1}^{\infty } \rbra{ 2^{k+1} - 1 } 2^{-k} \sum_{n=1}^{\infty } z^n \bP(\varphi^W_{-k}=n) 
= \sum_{k=1}^{\infty } \rbra{ 2^{k+1} - 1 } 2^{-k} (1-w)^k 
\label{} \\
=& \frac{2(1-w)}{w} - \frac{1-w}{1+w} 
\sim \sqrt{ \frac{2}{1-z} } 
\quad \text{as $ z \up 1 $}. 
\label{}
\end{align}
By the Tauberian theorem (see \cite[Proposition 4.2]{TZ}), we obtain 
\begin{align}
\sum_{n=1}^{N-1} c_n \sim \frac{\sqrt{2}}{\Gamma(3/2)}N^{1/2} 
= 2\sqrt{\frac{2}{\pi}}N^{1/2} 
\quad \text{as $ N \to \infty $}. 
\label{}
\end{align}

(i) 
By Lemma \ref{nthPF}, we have 
\begin{align}
w_N(\tilde{J}, \tilde{R}^{\pm}) 
=& \tilde{\mu}(\tilde{J} \cap \tilde{T}^{-1} \tilde{R}^{\pm}) 
+ \sum_{n=1}^{N-1} \int_{\tilde{X}} \tilde{T}^{\wedge n}_{\tilde{\mu}} 
1_{\tilde{J}_n \cap \tilde{R}^{\pm}}(\omega, x) \, \d \tilde{\mu} 
\label{} \\
=& \frac{1}{2} \mu \rbra{ \rbra{ \frac{1}{4},\frac{5}{8} } } 
+ \sum_{n=1}^{N-1} c_n \cdot \mu(I^{\pm}_0) 
= \frac{3}{16} + \frac{1}{4} \sum_{n=0}^{N-1} c_n 
\sim \frac{1}{2} \sqrt{\frac{2}{\pi}}N^{1/2} 
\label{}
\end{align}
as $ N \to \infty $. 
We also obtain 
\begin{align}
w_N(\tilde{J})
= w_N(\tilde{J}, \tilde{R}^+) 
+ w_N(\tilde{J}, \tilde{R}^-) 
+ \tilde{\mu}(\tilde{J} \cap \tilde{T}^{-1} \tilde{J}) 
\sim \sqrt{\frac{2}{\pi}}N^{1/2} 
\quad \text{as $ N \to \infty $}. 
\label{}
\end{align}

(iii) Again by Lemma \ref{nthPF}, we obtain 
\begin{align}
\frac{1}{w_N(\tilde{J}, \tilde{R}^{\pm})} 
\sum_{n=1}^{N-1} \tilde{T}^{\wedge n}_{\tilde{\mu}} 
1_{\tilde{J}_n \cap \tilde{R}^{\pm}}(\omega, x) 
= \frac{\sum_{n=1}^{N-1} c_n}{w_N(\tilde{J}, \tilde{R}^{\pm})} 
\cdot 1_{I^{\pm}_0}(x) \to 4 \cdot 1_{I^{\pm}_0}(x) 
\label{}
\end{align}
as $ N \to \infty $, where the convergence occurs uniformly in $ (\omega,x) $. 

(ii) This claim is obvious by (iii) and (i). 
}

Let us now complete the proof of Theorem \ref{arclaw}. 

\Proof[Proof of Theorem \ref{arclaw}]{
Proposition \ref{prop: Hata P} shows that 
the measure $ \mu $ given in \eqref{eq: im2} 
is a $ \lambda $-a.c. $ \sigma $-finite $ T $-invariant measure. 
Proposition \ref{thm: Hata ce} shows uniqueness. 

By Lemma \ref{lem: wandering}, we see that 
the assumptions of Theorems 3.1 and 3.2 of \cite{TZ} are all satisfied, 
and therefore we obtain the arcsine and Darling--Kac laws \eqref{eq: arclaw}. 
}

\section{Proof for the mBGI map} \label{mBGI}

The proof for the mBGI map is quite similar to that for the Hata map. 

\Proof[Proof of Theorem \ref{BGI DKlaw}]{
For the mBGI map with $ p=2/3 $, we have 
\begin{align}
\tau_1 I^-_k = 
\begin{cases}
I^-_{k-1} & (k \ge 1) \\
\sum_{j=0}^{\infty } I^-_j & (k=0) 
\end{cases}
, \quad & 
\tau_1 I^+_k = 
\begin{cases}
I^+_{k-1} & (k \ge 1) \\
\sum_{j=0}^{\infty } I^+_j & (k=0) 
\end{cases}
\label{} \\
\tau_2 I^-_k = I^+_{k+2} \ (k \ge 0) 
, \quad & 
\tau_2 I^+_k = I^-_{k+2} \ (k \ge 0) , 
\label{}
\end{align}
where $ I^-_k $'s and $ I^+_k $'s have been introduced in \eqref{eq: BGI alpha}. 
The transition probability $ q(s,t) $ defined in \eqref{eq: prmt} is given as 
\begin{align}
q(s,t) = 
\begin{cases}
\frac{1}{3} & \text{if $ (s,t) = (I^-_k,I^+_{k+2}) $ or $ (I^+_k,I^-_{k+2}) $ for $ k \ge 0 $} 
\\
\frac{2}{3} & \text{if $ (s,t) = (I^-_k,I^-_{k-1}) $ or $ (I^+_k,I^+_{k-1}) $ for $ k \ge 1 $} 
\\
\frac{2}{3} \cdot \frac{1}{2^{j+1}} & 
\text{if $ (s,t) = (I^-_0,I^-_j) $ or $ (I^+_0,I^+_j) $ for $ j \ge 0 $} 
\end{cases}
. 
\label{eq: BGI tp}
\end{align}
Noting that the random walk $ \{ W'_n \} $ on $ \bZ $ such that 
$ W'_n-W'_{n-1} = 2 $ with probability $ 1/3 $ 
and $ W'_n-W'_{n-1} = -1 $ with probability $ 2/3 $ 
is irreducible recurrent (see, e.g., \cite[Proposition 9.14]{Kal}), 
we can easily see that the transition matrix $ Q=(q(s,t))_{s,t \in \alpha } $ 
is irreducible recurrent. 
The $ Q $-invariant measure $ \rho $ on $ \alpha $ is characterized by the recurrence relations: 
\begin{align}
\rho(I^-_k) 
=& \frac{1}{3} \rho(I^+_{k-2}) + \frac{2}{3} \rho(I^-_{k+1}) + \frac{2}{3} \cdot \frac{1}{2^{k+1}} \rho(I^-_0) 
, \quad k=0,1,2,\ldots, 
\label{} \\
\rho(I^+_k) 
=& \frac{1}{3} \rho(I^-_{k-2}) + \frac{2}{3} \rho(I^+_{k+1}) + \frac{2}{3} \cdot \frac{1}{2^{k+1}} \rho(I^+_0) 
, \quad k=0,1,2,\ldots, 
\label{}
\end{align}
where we understand that $ \rho(I^{\pm}_{-1}) = \rho(I^{\pm}_{-2}) = 0 $. 
The measure $ \rho $ on $ \alpha $ defined by 
\begin{align}
\rho(I^-_k) = \rho(I^+_k) 
= \frac{3 \sqrt{2}}{8} 
\rbra{ \frac{8}{3} \cdot 2^k + \frac{1}{3} \cdot (-1)^k - 1 } \cdot 2^{-k-2} 
, \quad k=0,1,2,\ldots, 
\label{}
\end{align}
is easily proved to be a unique (up to a constant multiple) $ Q $-invariant measure. 

In a similar way to that for the Hata map, 
we can prove that the $ \tilde{\lambda} $-partition 
\begin{align}
\tilde{\alpha} 
= \cbra{ \tilde{s} : s \in \alpha } 
= \cbra{ \tilde{I}^-_k, \tilde{I}^+_k }_{k=0}^{\infty } , 
\label{}
\end{align}
is a Markov partition, and that 
the dynamical system $ (\tilde{X},\tilde{\mu},\tilde{T}) $ 
is conjugate to the irreducible recurrent Markov chain 
$ (\tilde{X},\tilde{\nu}_{\tilde{\mu}},\tilde{T}) $, 
where $ \tilde{\mu} $ is given as \eqref{eq: mmm}. 
Consequently, the dynamical system $ (\tilde{X},\tilde{\mu},\tilde{T}) $ 
is conservative ergodic, and so the measure $ \tilde{\mu} $ 
is a unique (up to a constant multiple) 
$ \tilde{\lambda} $-a.c. $ \sigma $-finite $ \tilde{T} $-invariant measure. 

For the junction, we may take 
\begin{align}
J := \rbra{ \frac{1}{4},\frac{3}{4} } =  I^-_0 + I^+_0 . 
\label{}
\end{align}
Then the Perron--Frobenius operator satisfies 
\begin{align}
\tilde{T}^{\wedge n}_{\tilde{\mu}} 1_{\tilde{J}_n}(\omega,x) = c'_n 1_J(x) 
, \quad \text{$ \tilde{\lambda} $-a.e.}, 
\label{eq: BGI tilde T wedge n} 
\end{align}
where 
\begin{align}
c'_n := \sum_{k=1}^{\infty } 
\frac{3 \sqrt{2}}{16} 
\rbra{ \frac{8}{3} + \frac{1}{3} \cdot \rbra{- \frac{1}{2}} ^k - \rbra{ \frac{1}{2} }^k } 
\bP(\varphi^{W'}_{-k} = n) . 
\label{eq: c'n}
\end{align}
From the formula \eqref{eq: BGI tilde T wedge n}, 
we can derive the wandering rate asymptotics 
\begin{align}
w_N(\tilde{J}) \sim \sqrt{\frac{2}{\pi}}N^{1/2} 
\quad \text{as $ N \to \infty $} , 
\label{eq: BGI wN tilde J}
\end{align}
which completes the proof by the help of Theorems 3.1 and 3.2 of \cite{TZ}. 

To obtain \eqref{eq: BGI wN tilde J}, 
we utilized the following formula 
instead of \eqref{eq: Hata w}: 
\begin{align}
w := 1 - \bE z^{\varphi^{W'}_{-1}} \sim \sqrt{1-z} 
\quad \text{as $ z \up 1 $} . 
\label{eq: mBGI w}
\end{align}
To obtain this formula, we note that 
\begin{align}
\bE z^{\varphi^{W'}_{-1}} 
= \frac{z}{3} \bE z^{\varphi^{W'}_{-3}} + \frac{2z}{3} 
= \frac{z}{3} \rbra{ \bE z^{\varphi^{W'}_{-1}} }^3 + \frac{2z}{3} 
, \quad 0<z<1 
\label{}
\end{align}
and hence $ w^2 - \frac{1}{3} w^3 = \frac{1-w}{z} \cdot (1-z) $, 
which yields 
\begin{align}
w = \sqrt{w^2} = \sqrt{ \frac{1}{1 - w/3} \cdot \frac{1-w}{z} \cdot (1-z) } 
\sim \sqrt{1-z} 
\quad \text{as $ z \up 1 $} , 
\label{}
\end{align}
since $ w \down 0 $ as $ z \up 1 $. 
}

\section{Proof for the Pelikan map} \label{Pelikan}

For the Pelikan map, the $ \tilde{\lambda} $-partition $ \tilde{\alpha } $ 
fails to be a Markov partition, and so we must take up another partition 
to reduce the problem to a Markov chain. 

\Proof[Proof of Theorem \ref{P DKlaw}]{
For the Pelikan map with $ p=1/2 $, we have 
\begin{align}
\tau_1 I_k = 
\begin{cases}
I_{k-1} \ (k \ge 1) 
\\
\sum_{j=0}^{\infty } I_j \ (k=0) 
\end{cases}
, \quad 
\tau_2 I_k = I_{k+1} \ (k \ge 0) , 
\label{}
\end{align}
where $ I_k $'s have been introduced in \eqref{eq: P alpha}. 
The transition probability $ q(s,t) $ defined in \eqref{eq: prmt} is given as 
\begin{align}
q(s,t) = 
\begin{cases}
\frac{1}{2} & \text{if $ (s,t) = (I_k,I_{k-1}) $ or $ (I_k,I_{k+1}) $ for $ k \ge 1 $} 
\\
\frac{1}{2^{j+2}} & \text{if $ (s,t) = (I_0,I_j) $ for $ j=0,2,3,\ldots $} 
\\
\frac{1}{2} + \frac{1}{2^{3}} & \text{if $ (s,t) = (I_0,I_1) $} 
\end{cases}
\label{}
\end{align}
We easily see that $ Q $ is irreducible 
and we can prove similarly to the Hata map that $ Q $ is recurrent. 
The $ Q $-invariant measure $ \rho $ on $ \alpha $ is characterized by 
\begin{align}
\rho(I_k) 
= \frac{1}{2} \rho(I_{k+1}) + \frac{1}{2} \rho(I_{k-1}) + \frac{1}{2^{k+2}} \rho(I_0) 
, \quad k=0,1,2,\ldots, 
\label{}
\end{align}
where we understand that $ \rho(I_{-1}) = 0 $, 
and the measure $ \rho $ on $ \alpha $ defined by 
\begin{align}
\rho(I_k) = 2^{-2} \cdot \rbra{ 2^{k+1} - 1 } \cdot 2^{-k-1} 
, \quad k=0,1,2,\ldots, 
\label{}
\end{align}
is easily proved to be a unique (up to a constant multiple) $ Q $-invariant measure. 

Unfortunately, the $ \tilde{\lambda} $-partition 
$ \tilde{\alpha } = \{ \tilde{s} : s \in \alpha \} $ 
is not a Markov partition for $ (\tilde{X},\tilde{\mu},\tilde{T}) $, 
because the map $ \tilde{T} : \tilde{I}_0 \to \tilde{T} \tilde{I}_0 $ is not invertible. 
To overcome this difficulty, we write 
\begin{align}
\tilde{s}^{\tau_1} = \{ \xi_1 = \tau_1 \} \cap \tilde{s} 
, \quad 
\tilde{s}^{\tau_2} = \{ \xi_1 = \tau_2 \} \cap \tilde{s} 
, \quad 
s \in \alpha . 
\label{}
\end{align}
By a similar argument to the Hata map using the recurrence of the simple symmetric random walk, 
we can prove that the $ \tilde{\lambda} $-partition of $ \tilde{X} $ defined by 
\begin{align}
\tilde{\alpha }^* = \cbra{ \tilde{I}_k^{\tau_1}, \ \tilde{I}_k^{\tau_2} }_{k=0}^{\infty } 
\label{}
\end{align}
is a Markov partition, and consequently the dynamical system 
$ (\tilde{X},\tilde{\mu},\tilde{T}) $ 
is conjugate to $ (\tilde{Y}^*,\tilde{\mu} \circ (\tilde{\phi}^*)^{-1},\theta) $, 
where $ \tilde{\mu} $ is given as \eqref{eq: mm}, 
and the graph shift $ \tilde{Y}^* $ is defined as 
\begin{align}
\tilde{Y}^* := \cbra{ \tilde{\vs}^* = (\tilde{s}^{\chi_0}_0,\tilde{s}^{\chi_1}_1,\tilde{s}^{\chi_2}_2,\ldots) : 
\tilde{s}^{\chi_0}_0,\tilde{s}^{\chi_1}_1,\tilde{s}^{\chi_2}_2,\ldots \in \tilde{\alpha}^* } 
\label{}
\end{align}
with the Borel bijection $ \tilde{\phi}^* : \tilde{X} \to \tilde{Y}^* $ defined as 
\begin{align}
\text{$ \tilde{\phi}^*(\omega,x) = \tilde{\vs}^* $ 
with $ \tilde{T}^n(\omega,x) \in \tilde{s}^{\chi_n}_n $ for all $ n $}. 
\label{}
\end{align}
The transition probability on $ \tilde{\alpha }^* $ defined as 
\begin{align}
\tilde{p}^*(\tilde{s}^{\chi_1},\tilde{t}^{\chi_2}) 
= \frac{\tilde{\lambda} \bigl( \, \tilde{s}^{\chi_1} \cap \tilde{T}^{-1} \tilde{t}^{\chi_2} \, \bigr) }
{\tilde{\lambda}(\tilde{s}^{\chi_1})} 
= \frac{1}{2} \frac{\lambda(s \cap \chi_1^{-1} t )}{\lambda(s)} 
, \quad s,t \in \alpha , \ \chi_1,\chi_2 \in \{ \tau_1, \tau_2 \} . 
\label{}
\end{align}
Noting that 
$ \tilde{p}^*(\tilde{s}^{\chi},\tilde{t}^{\tau_1}) 
= \tilde{p}^*(\tilde{s}^{\chi},\tilde{t}^{\tau_2}) $ and 
\begin{align}
\tilde{p}^*(\tilde{s}^{\tau_1},\tilde{t}^{\tau_1}) 
+ \tilde{p}^*(\tilde{s}^{\tau_2},\tilde{t}^{\tau_1}) 
= q(s,t) 
, \quad s,t \in \alpha , 
\label{}
\end{align}
we can prove that the transition matrix $ \tilde{P}^* $ is irreducible recurrent. 
We can then prove that the dynamical system $ (\tilde{X},\tilde{\mu},\tilde{T}) $ 
is conjugate to the irreducible recurrent Markov chain 
$ (\tilde{Y}^*,\tilde{\nu}^*_{\tilde{\mu}},\theta) $, 
where the measure $ \tilde{\nu}^*_{\tilde{\mu}} $ on $ \tilde{Y}^* $ is defined as 
\begin{align}
\tilde{\nu}^*_{\tilde{\mu}}([\tilde{s}^{\chi_0}_0,\tilde{s}^{\chi_1}_1,\ldots,\tilde{s}^{\chi_n}_n]) 
= \tilde{\mu}(\tilde{s}^{\chi_0}_0) 
\tilde{p}^*(\tilde{s}^{\chi_0}_0,\tilde{s}^{\chi_1}_1) \tilde{p}^*(\tilde{s}^{\chi_1}_1,\tilde{s}^{\chi_2}_2) 
\cdots \tilde{p}^*(\tilde{s}^{\chi_{n-1}}_{n-1},\tilde{s}^{\chi_n}_n) . 
\label{}
\end{align}

For the junction, we may take 
\begin{align}
J := \rbra{ \frac{1}{2},1 } = I_0 . 
\label{}
\end{align}
By a similar argument to that for the Hata map, we obtain 
\begin{align}
\tilde{T}^{\wedge n}_{\tilde{\mu}} 1_{\tilde{J}_n}(\omega,x) = 2 c_n 1_J(x) 
, \quad \text{$ \tilde{\lambda} $-a.e.}, 
\label{eq: P tilde T wedge n} 
\end{align}
where $ c_n $ is given in \eqref{eq: cn}, 
and obtain the wandering rate asymptotics 
\begin{align}
w_N(\tilde{J}) \sim \sqrt{\frac{2}{\pi}}N^{1/2} 
\quad \text{as $ N \to \infty $} , 
\label{eq: P wN tilde J}
\end{align}
which completes the proof by the help of Theorems 3.1 and 3.2 of \cite{TZ}. 
}

\end{document}